\documentclass[10pt]{amsart}
\usepackage{mathptmx}
\usepackage{amsmath}     
\usepackage{amssymb}
\usepackage{array}
\usepackage{geometry}
\usepackage[bookmarks=true,colorlinks=true, pdfstartview=FitV, linkcolor=black, citecolor=blue, urlcolor=black]{hyperref}
\usepackage{movie15}

\usepackage{color}
\definecolor{DarkRed}{rgb}{0.55,.00,0.2}
\definecolor{DarkGrey}{rgb}{0.35,.35,0.35}

\theoremstyle{definition}

\theoremstyle{remark}

\numberwithin{equation}{section}

\begin{document}

\title{A modification of the Prudnikov and Laguerre polynomials}

\author{{\bf S. Yakubovich} \\
 {\em {Department of Mathematics, Fac. Sciences of University of Porto,\\ Rua do Campo Alegre,  687; 4169-007 Porto (Portugal) }}}
  \thanks{ E-mail: syakubov@fc.up.pt} 
\vspace{3mm}

\thanks{ The work was  partially supported by CMUP, which is financed by national funds through FCT(Portugal),  under the project with reference UIDB/00144/2020. }

\subjclass[2000]{  33C47, 33C10, 44A15, 42C05 }

\date{\today}


\keywords{Orthogonal polynomials, modified Bessel function,  Mellin transform, associated Laguerre polynomials}

\begin{abstract} A two-parameter sequence of orthogonal polynomials $\{P_n( x; \lambda, t)\}_{n\ge 0}$  with respect to the weight function $x^\alpha e^{- \lambda x} \rho_\nu(x t),\  \alpha > -1,\  \lambda, t \ge 0, \ \rho_{\nu}(x)= 2 x^{\nu/2} K_\nu(2\sqrt x),\ x >0, \nu \ge 0$, where  $K_\nu(z)$ is the modified Bessel function, is investigated.  The case $\lambda=0$ corresponds to the Prudnikov polynomials and $t=0$ is related to the Laguerre polynomials. A special one-parameter case $\{P_n( x; 1-t, t)\}_{n\ge 0},\ t \in [0,1]$ is analyzed as well.  \\
\end{abstract}

\maketitle

\markboth{\rm \centerline{ S. Yakubovich}}{Orthogonal  polynomials  with non-classical weights}

\section{Introduction and preliminary results}

Recently \cite{YAP} the author gave an interpretation of the Prudnikov sequence of orthogonal polynomials  with  the weight function $x^\alpha \rho_\nu(x)$, where \  $ \rho_{\nu}(x)= 2 x^{\nu/2} K_\nu(2\sqrt x),\ x >0, \nu \ge 0,\ \alpha > -1$ and  $K_\nu(z)$ is  the modified Bessel function  \cite{Bateman}, Vol. II, in terms of the so-called composition orthogonality with respect to the measure $x^\nu e^{-x} dx$ related to the Laguerre polynomials \cite{Chi}.  Our main goal here  is to  characterize a two-parameter sequence of orthogonal polynomials $\{P_n( x; \lambda, t)\}_{n\ge 0}$  with respect to the weight function $x^\alpha e^{- \lambda x} \rho_\nu(x t),\  \alpha > -1,\  \lambda, t \ge 0$, satisfying the following orthogonality conditions

$$\int_{0}^{\infty}  P_{m} (x; \lambda, t) P_{n}(x; \lambda, t) e^{-\lambda x} \rho_\nu(x t) x^\alpha dx = \delta_{n,m},\eqno(1.1)$$
where $\delta_{n,m},\   n,m\in\mathbb{N}_{0}$ is the Kronecker symbol.     Up to a normalization factor conditions (1.1) are equivalent to the following equalities, respectively, 

$$\int_{0}^{\infty}  P_{n} (x; \lambda, t)  e^{-\lambda x} \rho_\nu(x t) \  x^{\alpha+m} dx = 0,\quad m=0,1,\dots, n-1,\ \quad n \in \mathbb{N}.\eqno(1.2)$$
The  function  $\rho_{\nu}$  has  the Mellin-Barnes integral representation
$$
\rho_\nu(x)=  \frac{1}{2\pi i} \int_{\gamma-i\infty}^{\gamma+i\infty} \Gamma(\nu+s) \Gamma (s) x^{-s} ds\ , \quad x, \gamma \in \mathbb{R}_{+},\ \nu  \in \mathbb{R},\eqno(1.3)
$$
where $\Gamma(z)$ is Euler's gamma-function \cite{Bateman}, Vol. I. Moreover, the Parseval equality for the Mellin transform \cite{Tit} permits to write (1.5) as the Laplace integral representation for this function. In fact, we obtain 

$$\rho_\nu(x)= \int_0^\infty y^{\nu-1} e^{-y - x/y} dy,\quad  x >0,\ \nu \in \mathbb{R}.\eqno(1.4)$$ 
Since $\rho_\nu(0) = \Gamma(\nu), \  \nu >0$, the polynomials $P_{n}(x; \lambda, 0),\  \lambda >0$ are associated with Laguerre polynomials, and, correspondingly,  $P_{n}(x; 0, t),\ t > 0,\ n \in \mathbb{N}_0$ are related to the Prudnikov polynomials \cite{YAP}. The asymptotic behavior of the modified Bessel function at infinity and near the origin \cite{Bateman}, Vol. II gives the corresponding values for the  function $\rho_\nu,\ \nu \in \mathbb{R}$.  Precisely, we have
$$\rho_\nu (x)= O\left( x^{(\nu-|\nu|)/2}\right),\  x \to 0,\ \nu\neq 0, \quad  \rho_0(x)= O( \log x),\ x \to 0,\eqno(1.5)$$

$$ \rho_\nu(x)= O\left( x^{\nu/2- 1/4} e^{- 2\sqrt x} \right),\ x \to +\infty.\eqno(1.6)$$
The moments for the weight $e^{- \lambda x} \rho_\nu(xt),\ \lambda^2+ t^2 \neq 0$ can be calculated via  integral  (1.4).  Indeed, we have 

$$I_n= \int_0^\infty e^{-\lambda x}  \rho_\nu(x t) x^{n +\alpha} dx =  \lambda^{-n -\alpha-1} \int_0^\infty e^{-x}   x^{n+\alpha}  \int_0^\infty u^{\nu-1} e^{-u - x t/(\lambda u)} du dx $$

$$=  \Gamma(1+n+\alpha)    \int_0^\infty {u^{\nu+n+\alpha}  e^{-u } \over (\lambda u +t)^{n+\alpha+1} }  du,\quad  n \in \mathbb{N}_0.\eqno(1.7)$$
where the interchange of the order of integration is justified by Fubini's theorem.  The latter integral  is expressed in terms of the Tricomi function $\Psi(a,b; z)$ (cf. \cite{Bateman},  Vol. I and \cite{PrudnikovMarichev}, Vol. I, Entry 2.3.6.9).  Hence we get finally the values from (1.7) 

$$I_n= \lambda^{-\nu-\alpha-n-1} t^\nu \Gamma(n+\nu+\alpha+1)\Gamma(n+\alpha+1) \Psi \left(1+n+\alpha+\nu,\ 1+\nu;  {t\over \lambda} \right).\eqno(1.8)$$
The quotient of the scaled Macdonald functions $\rho_\nu, \rho_{\nu+1}$ is given by the Ismail integral \cite{Ismail}

$${\rho_\nu(x) \over \rho_{\nu+1}(x) } = {1\over \pi^2} \int_0^\infty {y^{-1} dy \over (x+y) \left[ J_{\nu+1}^2 (2\sqrt y)+ Y^2_{\nu+1}(2\sqrt y) \right] },\eqno(1.9)$$
where $J_\nu(z), Y_\nu(z)$ are Bessel functions of the first and second kind, respectively \cite{Bateman}.  Moreover, Entry 2.19.4.13  in \cite{PrudnikovMarichev}, Vol. II represents the product $x^n \rho_\nu(x)$ in terms of the associated Laguerre polynomials \cite{Chi}.  Precisely, it has

$${(-1)^n x^n\over n!}\  \rho_\nu(x)=   \int_0^\infty y^{\nu+n -1} e^{-y - x/y}  L_n^\nu(y) dy,\quad    n \in\mathbb{N}_{0}.\eqno(1.10)$$
Appealing to   the Riemann-Liouville fractional integral \cite{YaL}

$$ \left( I_{-}^\alpha  f \right) (x)  = {1\over \Gamma(\alpha)} \int_x^\infty (t-x)^{\alpha-1} f(t) dt,\quad  {\rm Re} \alpha > 0\eqno(1.11)$$
and Entry 2.16.3.8  in \cite{Bateman}, Vol. II, we get  the  formula 
$$\rho_\nu(x)= \left( I_{-}^\nu \rho_0 \right) (x).\eqno(1.12)$$
Besides, the index law for fractional integrals immediately implies

$$ \rho_{\nu+\mu} (x)= \left( I_{-}^\nu \rho_\mu \right) (x)=   \left( I_{-}^\mu \rho_\nu \right) (x).\eqno(1.13)$$
The corresponding definition of the fractional derivative presumes the relation $ D^\mu_{-}= - D  I_{-}^{1-\mu}$.   Hence for the ordinary $n$-th derivative of $\rho_\nu$ we find

$$D^n \rho_\nu(x)= (-1)^n \rho_{\nu-n} (x),\quad n \in \mathbb{N}_0.\eqno(1.14)$$
 Recalling  the Mellin-Barnes integral (1.3) and reduction formula for the gamma-function it is not difficult to derive the recurrence relation 
 
 $$\rho_{\nu+1} (x) =    \nu \rho_\nu(x)+ x \rho_{\nu-1} (x),\quad \nu \in \mathbb{R}.\eqno(1.15)$$
As it follows from the theory of orthogonal polynomials \cite{Chi},  a sequence $\{p_n\}_{n\ge 0}$ satisfies the three term recurrence relation 

$$x p_n(x) = A_{n+1} p_{n+1}(x) +B_n p_n(x) + A_{n} p_{n-1}(x),\eqno(1.16)$$
where $p_{-1}(x) \equiv 0,\  p_n(x) = a_n x^n+ b_n  x^{n-1}+ \dots,$ and 

$$A_n= {a_{n-1} \over a_n},\quad B_n= {b_{n}\over a_n} - {b_{n+1}\over a_{n+1}}.\eqno(1.17)$$
As a consequence of (1.16) the Christoffel-Darboux formula takes place

$$\sum_{k=0}^n   p_k(x) p_k(y) = A_{n+1} \frac{  p_{n+1}(x) p_n(y) -  p_{n}(x) p_{n+1}(y)}{x-y}.\eqno(1.18)$$
Denoting the weight function as $\omega(x)=e^{-\lambda x} \rho_\nu(x t) x^\alpha$ we establish a differential equation whose solution is $\omega$.

{\bf Lemma 1}. {\it The weight function $ \omega$ satisfies the following second order differential equation} 

$$x^2\omega^{\prime\prime}-  (2(\alpha- \lambda x)  +\nu-1) x \omega^\prime + \left(  \left( \alpha -\lambda x\right)^2 + x( \lambda(1-\nu) -1) + \alpha \nu \right) \omega =0.\eqno(1.19)$$

\begin{proof} In fact, recalling (1.14), (1.15), we have

$$[xt  \omega ]^\prime =  {d\over dx} \left[ x^\alpha e^{- \lambda x} (\rho_{\nu+2} (xt)- (\nu+1)\rho_{\nu+1}(xt)) \right]$$

$$=t \left(  \alpha -\lambda x  + \nu+1\right) \omega - e^{-\lambda x} \rho_{\nu+1} (x t) x^\alpha.\eqno(1.20)$$
Hence with one more differentiation it gives

$$[x t  \omega ]^{\prime\prime} = t(1- \lambda) \omega +  t \left(  \alpha -\lambda x  + \nu+1\right) \omega^\prime -
\left(  {\alpha\over x} -\lambda \right) e^{-\lambda x} \rho_{\nu+1} (x t) x^\alpha.\eqno(1.21) $$
Combining (1.20), (1.21), we find after simplification 

$$ [x  \omega ]^{\prime\prime} =   \left(  2( \alpha -\lambda x)   + \nu+1\right) \omega^\prime - \left(  \left(  {\alpha\over x} -\lambda \right) \left(   \alpha -\lambda x  + \nu\right) +\lambda -1 \right) \omega.$$
Therefore fulfilling the differentiation and multiplying by $x$, we arrive at (1.19). 

\end{proof}

\section{Differential-difference properties}

We begin this section with the following result.

{\bf Theorem 1}. {\it Let $\alpha > 0, \ \lambda, t, \nu  \ge 0, \ \lambda^2+t^2 \neq 0, n \in \mathbb{N}$. Orthogonal polynomials $P_n (x; \lambda, t)$ satisfy the integro-differential-difference  equation

$$x {\partial\over \partial x} [P_n (x; \lambda, t)]  =  = { x t A_n\over \pi^2 }  \int_0^\infty   \int_0^\infty   \bigg[P_n (x; \lambda,t)  P_{n-1} (y; \lambda, t) -  P_n (y; \lambda,t)  P_{n-1} (x; \lambda, t) \bigg] $$

$$\times    {  P_n (y; \lambda,t) \ e^{-\lambda y} \rho_\nu(y t) y^{\alpha}  \over u  (y+u) (x+u) \left[ J_{\nu}^2 (2\sqrt u)+ Y^2_{\nu}(2\sqrt u) \right] } dudy $$

$$ -  \alpha A_n  \bigg[a_{n-1,0}\  P_n (x; \lambda,t)   -  a_{n,0} \ P_{n-1} (x; \lambda, t) \bigg] \int_0^\infty  P_n (y; \lambda,t)   e^{-\lambda y} \rho_\nu(y t) y^{\alpha-1} dy,\eqno(2.1)$$
where $A_n$ is defined in $(1.17),  P_n (x; \lambda, t) = a_n x^n+ b_n  x^{n-1}+ \dots + a_{n,0}$ and, clearly,  polynomial coefficients are functions of $\lambda, t$}.

\begin{proof} Since $ {\partial\over \partial x} [P_n (x; \lambda, t)] $ is a polynomial of degree $n-1$, we write it in the form

$${\partial\over \partial x} [P_n (x; \lambda, t)]  = \sum_{k=0}^{n-1} c_{n,k}  P_k (x; \lambda, t),\eqno(2.2)$$
where, owing to the orthogonality,

$$c_{n,k} = \int_0^\infty  {\partial\over \partial x} [P_n (x; \lambda, t)] P_k (x; \lambda, t) \ e^{-\lambda x} \rho_\nu(x t) x^\alpha dx.\eqno(2.3)$$
Then, integrating by parts with the use of  the orthogonality (1.1) and relation (1.14), we eliminate integrated terms by virtue of the asymptotic formulas (1.5), (1.6) to obtain

$$c_{n,k} = t \int_0^\infty  P_n (x; \lambda,t) P_k (x; \lambda, t) \  e^{-\lambda x} \rho_{\nu-1} (x t) x^\alpha dx $$

$$- \alpha \int_0^\infty  P_n (x; \lambda,t) P_k (x; \lambda, t) \  e^{-\lambda x} \rho_\nu(x t) x^{\alpha-1} dx.\eqno(2.4)$$
In the meantime, we easily observe from the orthogonality that 

$$ \int_0^\infty  P_n (y; \lambda, t) \sum_{k=0}^{n-1}  P_k (y; \lambda, t) P_k (x; \lambda, t)   e^{-\lambda y} \rho_\nu(y t) y^{\alpha} \  {\rho_{\nu-1} (x t) \over \rho_{\nu} (x t)  }dy = 0,$$

$$ \int_0^\infty  P_n (y; \lambda, t) \sum_{k=0}^{n-1}  P_k (y; \lambda, t) P_k (x; \lambda, t)   e^{-\lambda y} \rho_\nu(y t) y^{\alpha}  x^{-1} dy = 0.$$

Therefore from (1.2), (1.9), (2.2), (2.4) and Christoffel-Darboux formula (1.18) we derive

$${\partial\over \partial x} [P_n (x; \lambda, t)]  = -\alpha  \sum_{k=0}^{n-1} P_k (x; \lambda, t)  \int_0^\infty  P_n (y; \lambda,t) P_k (y; \lambda, t) \  e^{-\lambda y} \rho_\nu(y t) y^\alpha \left[{1\over y} - {1\over x} \right] dy$$

$$+t  \sum_{k=0}^{n-1} P_k^\nu(x,t) \int_0^\infty   P_n (y; \lambda,t) P_k (y; \lambda, t) \  e^{-\lambda y} \rho_\nu(y t) y^\alpha  \bigg[  {\rho_{\nu-1} (y t) \over \rho_{\nu} (y t)  } -  {\rho_{\nu-1} (x t) \over \rho_{\nu} (x t)  } \bigg] dy$$

$$= -  {\alpha A_n\over x} \   \int_0^\infty  P_n (y; \lambda,t) \bigg[P_n (x; \lambda,t)  P_{n-1} (y; \lambda, t) -  P_n (y; \lambda,t)  P_{n-1} (x; \lambda, t) \bigg] \  e^{-\lambda y} \rho_\nu(y t) y^{\alpha-1} dy  $$

$$ +  t A_n \int_0^\infty   P_n (y; \lambda,t) \bigg[P_n (x; \lambda,t)  P_{n-1} (y; \lambda, t) -  P_n (y; \lambda,t)  P_{n-1} (x; \lambda, t) \bigg] $$

$$\times  {e^{-\lambda y} \rho_\nu(y t) y^{\alpha}\over x-y}\bigg[  {\rho_{\nu-1} (y t) \over \rho_{\nu} (y t)  } -  {\rho_{\nu-1} (x t) \over \rho_{\nu} (x t)  } \bigg]  dy  $$

$$ = {t A_n\over \pi^2 }  \int_0^\infty   \int_0^\infty   \bigg[P_n (x; \lambda,t)  P_{n-1} (y; \lambda, t) -  P_n (y; \lambda,t)  P_{n-1} (x; \lambda, t) \bigg] $$

$$\times    {  P_n (y; \lambda,t) \ e^{-\lambda y} \rho_\nu(y t) y^{\alpha}  \over u  (y+u) (x+u) \left[ J_{\nu}^2 (2\sqrt u)+ Y^2_{\nu}(2\sqrt u) \right] } dudy $$

$$ -  {\alpha A_n\over x} \  \bigg[a_{n-1,0}\  P_n (x; \lambda,t)   -  a_{n,0} \ P_{n-1} (x; \lambda, t) \bigg] \int_0^\infty  P_n (y; \lambda,t)   e^{-\lambda y} \rho_\nu(y t) y^{\alpha-1} dy.$$
This completes the proof of Theorem 1. 

\end{proof} 

Now, assuming that polynomial coefficients are continuously differentiable functions of $\lambda, t$, we differentiate through equalities (1.2) with respect to $\lambda$ and $t$ under the integral sign and integrate by parts with the use of the orthogonality and equalities (1.14), (1.15). We note that the differentiation under the integral sign in (1.2) can be motivated by the absolute and uniform convergence in the subset of $\mathbb{R}^2\   M= \left\{ (\lambda, t ) \in \mathbb{R}^2_+ \left| \right. \lambda \ge \lambda_0 > 0,\  t \ge t_0 > 0\right\}.$   Hence we obtain from (1.2) three types of orthogonality conditions

$$\int_{0}^{\infty}  {\partial\over \partial \lambda } [ P_{n} (x; \lambda, t) ] e^{-\lambda x} \rho_\nu(x t) \  x^{\alpha+m} dx$$

$$ -  \int_{0}^{\infty}   P_{n} (x; \lambda, t) e^{-\lambda x} \rho_\nu(x t) \  x^{\alpha+m+1} dx = 0,\quad m=0,1,\dots, n-1,\eqno(2.5)$$

$$\int_{0}^{\infty}  t {\partial\over \partial t} [ P_{n} (x; \lambda, t) ] e^{-\lambda x} \rho_\nu(x t) \  x^{\alpha+m} dx $$

$$-  \int_{0}^{\infty}   P_{n} (x; \lambda, t) e^{-\lambda x} \rho_{\nu+1}(x t) \  x^{\alpha+m} dx = 0,\quad m=0,1,\dots, n-1,\eqno(2.6)$$

$$ \lambda \int_{0}^{\infty}   P_{n} (x; \lambda, t) e^{-\lambda x} \rho_\nu(x t) \  x^{\alpha+m+1} dx +  \int_{0}^{\infty}   P_{n} (x; \lambda, t) e^{-\lambda x} \rho_{\nu+1} (x t) \  x^{\alpha+m} dx $$

$$- \int_{0}^{\infty}  {\partial\over \partial x } [ P_{n} (x; \lambda, t) ] e^{-\lambda x} \rho_\nu(x t) \  x^{\alpha+m+1} dx = 0,\quad m=0,1,\dots, n-1.\eqno(2.7)$$
Then as a direct consequence of (2.6), (2.7) we find the equalities

$$ \int_{0}^{\infty}  \left( t {\partial\over \partial t} -  x {\partial\over \partial x } \right)  [ P_{n} (x; \lambda, t) ] e^{-\lambda x} \rho_\nu(x t) \  x^{\alpha+m} dx$$

$$+  \lambda \int_{0}^{\infty}   P_{n} (x; \lambda, t) e^{-\lambda x} \rho_\nu(x t) \  x^{\alpha+m+1} dx = 0,\quad m=0,1,\dots, n-1.\eqno(2.8)$$
Furthermore,  (1.16),  (2.5), (2.8) and the differentiation through (1.1) for $n=m$ by $\lambda$ and by $t$  yield 

$$\int_{0}^{\infty} \bigg[  {\partial\over \partial \lambda } [ P_{n} (x; \lambda, t) ] -  A_n P_{n-1} (x; \lambda, t) \bigg] e^{-\lambda x} \rho_\nu(x t) \  x^{\alpha+m} dx= 0,\quad m=0,1,\dots, n-1,\eqno(2.9)$$

$$ \int_{0}^{\infty}  \bigg[ \left( t {\partial\over \partial t} -  x {\partial\over \partial x } \right)  [ P_{n} (x; \lambda, t) ]  +\lambda A_n P_{n-1} (x; \lambda, t) \bigg] e^{-\lambda x} \rho_\nu(x t) \  x^{\alpha+m} dx =  0,\quad m=0,1,\dots, n-1,\eqno(2.10)$$

$$2 \int_{0}^{\infty}  P_{n} (x; \lambda, t)  {\partial\over \partial \lambda } [ P_{n} (x; \lambda, t) ] e^{-\lambda x} \rho_\nu(x t) \  x^{\alpha} dx-  \int_{0}^{\infty}   [ P_{n} (x; \lambda, t) ]^2  e^{-\lambda x} \rho_\nu(x t) \  x^{\alpha+1} dx = 0,\eqno(2.11)$$

$$2 \int_{0}^{\infty}  P_{n} (x; \lambda, t)  {\partial\over \partial t } [ P_{n} (x; \lambda, t) ] e^{-\lambda x} \rho_\nu(x t) \  x^{\alpha} dx-  \int_{0}^{\infty}   [ P_{n} (x; \lambda, t) ]^2  e^{-\lambda x} \rho_{\nu-1} (x t) \  x^{\alpha+1} dx = 0.\eqno(2.12)$$
As a result we are ready to prove the following theorems.

{\bf Theorem 2.} {\it Let $\alpha > -1,\ \lambda, t, \nu  \ge 0, \ \lambda^2+t^2 \neq 0, n \in \mathbb{N}$.  Orthogonal polynomials $P_{n} (x; \lambda, t) $ satisfy the  first order ordinary and partial differential-difference equations}

$$ {\partial\over \partial \lambda } [ P_{n} (x; \lambda, t) ] - {1\over a_n} {\partial a_n \over \partial \lambda } P_{n} (x; \lambda, t) -  A_n P_{n-1} (x; \lambda, t) = 0,\eqno(2.13)$$

$$ \left( t {\partial\over \partial t} -  x {\partial\over \partial x } \right)  [ P_{n} (x; \lambda, t) ]  - \bigg[ { t \over a_n}\  {\partial a_n \over \partial t } -  n \bigg] P_{n} (x; \lambda, t) +\lambda A_n P_{n-1} (x; \lambda, t) = 0.\eqno(2.14)$$

\begin{proof}    Indeed, equations (2.13), (2.14) are direct consequences of the uniqueness and  orthogonality conditions (2.9), (2.10), respectively. In fact, polynomials within the brackets under the integral sign are of degree $n$ and up to  constants are equal to $P_n$.  These constants, in turn,  can be found, equating coefficients  in front of $x^n$. 

\end{proof}
Equating coefficients in front of $x^{n-1}$ in (2.13), (2.14) to zero and dividing the obtained equalities by $a_n$, we find after simplification 

{\bf Corollary 1.}  {\it The following differential recurrence relations for coefficients $a_n, b_n, A_n, B_n$   hold}

$$   {\partial  \over \partial \lambda } \left( {b_n \over a_n} \right)  = A_n^2,\eqno(2.15)$$

$$   {\partial  \over \partial t } \left(  {t b_n \over a_n} \right)  = - \lambda A_n^2,\eqno(2.16)$$

$$\left( \lambda {\partial  \over \partial \lambda } +  t {\partial  \over \partial t } \right) B_n +  B_n = 0,\eqno(2.17)$$
$$\left( \lambda {\partial  \over \partial \lambda } -  t {\partial  \over \partial t } \right) \left( {b_n \over a_n} \right) -  {b_n \over a_n} = 2 \lambda A_n^2,\eqno(2.18)$$

$$\left( \lambda {\partial  \over \partial \lambda } -  t {\partial  \over \partial t } \right) B_n -  B_n = 2 \lambda \left[ A_n^2- A_{n+1}^2\right].\eqno(2.19)$$

We note that relations (2.17)-(2.19) follow immediately from (2.15), (2.16) via (1.17).  Moreover, as we will see,  identities (2.11), (2.12) generate 

{\bf Theorem 3}. {\it Let $\alpha > -1,\  t >0,\  \lambda, \nu  \ge 0, \  n \in \mathbb{N}_0$.  Then it has the equalities}

$$B_n =  {2\over a_n} {\partial a_n \over \partial \lambda },\eqno(2.20)$$

$$  {2 t\over a_n}  {\partial a_n \over \partial t} =  2n + \alpha+1 - \lambda B_n,\eqno(2.21)$$

$$\left( \lambda {\partial  \over \partial \lambda } +  t {\partial  \over \partial t } \right) a_n  = \left( n + {\alpha +1\over 2}\right) a_n,\eqno(2.22)$$

$$\left( \lambda {\partial  \over \partial \lambda } +  t {\partial  \over \partial t } \right) A_n  + A_n = 0.\eqno(2.23)$$

\begin{proof}  Identity (2.20) follows from (2.11),  orthogonality (1.1), the three term recurrence relation (1.16) and, as a consequence, from  the integral

$$\int_{0}^{\infty}  P_{n} (x; \lambda, t)  e^{-\lambda x} \rho_\nu(x t) \  x^{\alpha+n} dx = {1\over a_n}.\eqno(2.24)$$
In order to prove (2.21), we have from (2.12) via (1.14), (1.15), (2.24) and integration by parts

$$ {2\over a_n}  {\partial a_n \over \partial t} - {\alpha+1\over t} + {\lambda\over t} \ B_n - {2\over t}  \int_{0}^{\infty}  P_{n} (x; \lambda, t)  {\partial\over \partial x} [ P_{n} (x; \lambda, t) ] e^{-\lambda x} \rho_\nu(x t) \  x^{\alpha+1} dx = 0.\eqno(2.25)$$
Using again (2.24),   it is easily seen that the latter integral is equal to $n$.  Hence   we end up with  (2.21), and equalities  (2.22) , (2.23)  are  immediate consequences of (2.20), (2.21).

\end{proof}

Further, differentiating through by $\lambda$ and by $t$ the three term recurrence relation (1.16) and employing (2.17), (2.23), we arrive at the following result.

{\bf Corollary 2.} {\it Let $\alpha > -1,\ \lambda, t, \nu  \ge 0, \ \lambda^2+t^2 \neq 0,   n \in \mathbb{N}_0, m \in \mathbb{Z}$.  Orthogonal polynomials $P_n(x; \lambda, t)$  obey the partial differential-difference recurrence relation of the form

$$ \left(  \lambda {\partial  \over \partial \lambda } +  t {\partial  \over \partial t }  +  m \mathbb{E} \right) x P_n(x; \lambda, t)  = A_{n+1} \left( \lambda {\partial  \over \partial \lambda } +  t {\partial  \over \partial t } + (m-1)  \mathbb{E}\right) P_{n+1} (x; \lambda, t) $$

$$+B_n \left( \lambda {\partial  \over \partial \lambda } +  t {\partial  \over \partial t } + (m-1) \mathbb{E} \right) P_n(x; \lambda, t) + A_{n} \left( \lambda {\partial  \over \partial \lambda } +  t {\partial  \over \partial t } + (m-1) \mathbb{E} \right) P_{n-1}(x; \lambda, t),\eqno(2.26)$$
where $\mathbb{E}$ is the identity operator.}

The differential-difference recurrence relation (2.26) can be extended for the integer power of the operator $  \lambda {\partial  \over \partial \lambda } +  t {\partial  \over \partial t }  +  m \mathbb{E}, \  m \in \mathbb{Z}$.  Namely,  employing (2.17), (2.23) and the method of mathematical induction, it is not difficult to establish  the following    differential-difference recurrence relation 

$$ \left(  \lambda {\partial  \over \partial \lambda } +  t {\partial  \over \partial t }  +  m \mathbb{E} \right)^r x P_n(x; \lambda, t)  = A_{n+1} \left( \lambda {\partial  \over \partial \lambda } +  t {\partial  \over \partial t } + (m-1) \mathbb{E}\right)^r P_{n+1} (x; \lambda, t) $$

$$+B_n \left( \lambda {\partial  \over \partial \lambda } +  t {\partial  \over \partial t }+ (m-1) \mathbb{E} \right)^r P_n(x; \lambda, t) + A_{n} \left( \lambda {\partial  \over \partial \lambda } +  t {\partial  \over \partial t } + (m-1) \mathbb{E}\right)^r P_{n-1}(x; \lambda, t),\eqno(2.27)$$
where $r \in \mathbb{N}_0$.

Now, recalling Theorem 2, we multiply (2.13) by $\lambda$ and sum with (2.14) to get

{\bf Corollary 3.} {\it Let $\alpha > -1,\ \lambda, t, \nu > 0, \ n \in \mathbb{N}_0$.  Orthogonal polynomials $P_n(x; \lambda, t)$  satisfy the partial differential equation of the first order

$$ \left(  \lambda {\partial  \over \partial \lambda } +  t {\partial  \over \partial t }   \right)  P_n(x; \lambda, t)  - x  {\partial\over \partial x }   [ P_{n} (x; \lambda, t) ] -  {\alpha +1\over 2}  P_n(x; \lambda, t) = 0.\eqno(2.28)$$
Moreover,  analogously to $(2.27)$ for the operator powers it has the relation}

$$ \left(  \lambda {\partial  \over \partial \lambda } +  t {\partial  \over \partial t }  -   {\alpha +1\over 2} \right)^m  P_n(x; \lambda, t) = \left(x  {\partial\over \partial x }\right)^m     P_{n} (x; \lambda, t), \quad m \in \mathbb{N}_0.\eqno(2.29)$$

Equating coefficient in front of $x^{n-1}$ to zero, we find a companion of the equality (2.22)

$$\left( \lambda {\partial  \over \partial \lambda } +  t {\partial  \over \partial t } \right) b_n  = \left( n + {\alpha -1\over 2}\right) b_n.\eqno(2.30)$$
Further use of the three term recurrence relation (1.16) and the orthogonality (1.1) generates more values of the integrals similar to (2.24) and identities between polynomial coefficients.  In fact, we have

{\bf Lemma 2}. {\it Writing the polynomial $P_n$ in the form

$$ P_{n} (x; \lambda, t) = a_n x^n + b_n x^{n-1} + d_n x^{n-2} + \hbox{lower degrees},$$
the following values of integrals take place

$$ \int_{0}^{\infty}  [ P_{n} (x; \lambda, t) ]^2 e^{-\lambda x} \rho_\nu(x t) \  x^{\alpha+2} dx = A_{n+1}^2+ B_n^2 + A_n^2,\eqno(2.31)$$

$$\int_{0}^{\infty}  P_{n} (x; \lambda, t)  e^{-\lambda x} \rho_\nu(x t) \  x^{\alpha+n+1} dx = - {b_{n+1} \over a_{n+1} a_n},\eqno(2.32)$$

$$\int_{0}^{\infty}  P_{n} (x; \lambda, t)  e^{-\lambda x} \rho_\nu(x t) \  x^{\alpha+n+2} dx = {b_{n+2} b_{n+1}\over a_{n+2} a_{n+1} a_n}- {d_{n+2} 
\over a_{n+2}  a_n},\eqno(2.33)$$

 $$\int_{0}^{\infty}   [ P_{n} (x; \lambda, t) ]^2  e^{-\lambda x} \rho_\nu(x t) \  x^{\alpha+3} dx = A_{n+1}^2 \left[ B_{n+1} + 2 B_n\right] + A_n^2  B_{n-1}  + \left[ 2A_n^2 + B^2_n \right] B_n\eqno(2.34)$$
and the identity}
$$ {d_{n} \over a_{n} } - {d_{n+2} \over a_{n+2} } - {b_{n+1}\over a_{n+1}} \bigg[ B_n+ B_{n+1}\bigg]  = A_{n+1}^2+ B_n^2 + A_n^2.\eqno(2.35)$$

\begin{proof} Identity (2.31) is a direct consequence of (1.16), (1.1).   In the same manner one gets (2.32), (2.33), involving the polynomial $P_{n+1},\ P_{n+2}$, respectively. Then a combination with (2.31) and (1.17) yields (2.35). Finally,    the value of the integral (2.34) can be obtained, writing it via (1.16), namely,

$$ \int_{0}^{\infty}   [ P_{n} (x; \lambda, t) ]^2  e^{-\lambda x} \rho_\nu(x t) \  x^{\alpha+3} dx $$

$$=   \int_{0}^{\infty}   [ A_{n+1} P_{n+1} (x; \lambda, t) +B_n P_n(x; \lambda, t) + A_{n} P_{n-1}(x; \lambda, t) ]^2  e^{-\lambda x} \rho_\nu(x t) \  x^{\alpha+1} dx, $$
and appealing to (1.1), (1.16).

\end{proof}

On the other hand, we observe from (2.12), (2.21)

$$ \int_{0}^{\infty}   [ P_{n} (x; \lambda, t) ]^2  e^{-\lambda x} \rho_{\nu+1} (x t) \  x^{\alpha} dx = 2n+\alpha+\nu+1- \lambda B_n.\eqno(2.36)$$
Then, we employ (1.14), (1.15), (2.31), (2.32)  and integration by parts to get

$$ \int_{0}^{\infty}   [ P_{n} (x; \lambda, t) ]^2  e^{-\lambda x} \rho_{\nu+1} (x t) \  x^{\alpha+1} dx = \nu B_n +  t \int_{0}^{\infty}   [ P_{n} (x; \lambda, t) ]^2  e^{-\lambda x} \rho_{\nu-1} (x t) \  x^{\alpha+2} dx$$

$$= \left( \alpha +2 + \nu \right) B_n - \lambda \left[ A_{n+1}^2+ B_n^2 + A_n^2\right] + 2 \int_{0}^{\infty}    P_{n} (x; \lambda, t)  {\partial\over \partial x }   [ P_{n} (x; \lambda, t) ]   e^{-\lambda x} \rho_{\nu} (x t) \  x^{\alpha+2} dx$$

$$=   \left( \alpha +2 + \nu \right) B_n - \lambda \left[ A_{n+1}^2+ B_n^2 + A_n^2\right] + 2(n-1) {b_n\over a_n} - 2 n {b_{n+1} \over a_{n+1}},$$
i.e. we derive the identity

$$ \int_{0}^{\infty}   [ P_{n} (x; \lambda, t) ]^2  e^{-\lambda x} \rho_{\nu+1} (x t) \  x^{\alpha+1} dx =  \left( \alpha +2(n+1)  + \nu \right) B_n - \lambda \left[ A_{n+1}^2+ B_n^2 + A_n^2\right] - 2  {b_{n} \over a_{n}}.\eqno(2.37)$$

Now, returning to equalities (2.8), we observe that for $\lambda \neq 0$ the second integral is zero when $m=0,1,\dots, n-2$.  Consequently, if $t=0$,  we end up with a modification of Laguerre polynomials $\hat{L}^{\alpha}_n$, getting the known property

$${d \over dx }   [ P_{n} (x; \lambda, 0) ] = \hat{L}^{\alpha+1}_{n-1} (x; \lambda),\eqno(2.38)$$
and can be calculated explicitly, invoking properties of Laguerre polynomials.  Otherwise, in the case $t\neq 0$,    we write (2.8) in the form 

$$ \int_{0}^{\infty}  \left( t {\partial\over \partial t} -  x {\partial\over \partial x } \right)  [ P^\alpha_{n} (x; \lambda, t) ] e^{-\lambda x} \rho_\nu(x t) \  x^{\alpha+m} dx = 0,\quad m=0,1,\dots, n-2,\eqno(2.39)$$
$$ \int_{0}^{\infty}  \left( t {\partial\over \partial t} -  x {\partial\over \partial x } \right)  [ P^\alpha_{n} (x; \lambda, t) ] e^{-\lambda x} \rho_\nu(x t) \  x^{\alpha+n-1} dx \neq  0.\eqno(2.40)$$
This means  that the sequence $\left\{\left( t {\partial\over \partial t} -  x {\partial\over \partial x } \right)  [ P_{n} (x; \lambda, t) ] \right\}_{n\ge 0}$ is quasi-orthogonal with respect to the weight $e^{-\lambda x} \rho_\nu(x t) x^{\alpha}$. The latter property is confirmed, in particular, by virtue of the recurrence relation (2.14) as a combination of the above quasi-orthogonal polynomials of degree $n$ in terms of  their orthogonal counterparts  $P_n, P_{n-1}$. 

\section{Integral-diference equations. Composition orthogonality}

In this section we will establish the integral-difference equations for the sequence $\{ P_{n} \}_{n\ge 0}$ and characterize it in terms of the composition orthogonality for the operator polynomials $P_n(\theta)$, where $\theta = x  D x,\ D= {d\over dx} $ and the usual product of polynomials $P_n(x)P_m(x)$ is substituted by the composition of operators $P_n(\theta)P_m(\theta)$.  This notion was introduced in \cite{YAP} and as we will see below,  the orthogonality relations (1.2) can be rewritten in the composition sense by virtue of the key property for the operator $\theta$, namely, its nonnegative integer power can be represented  in the form 

$$\theta^m = x^m  D^m x^m,\quad m \in \mathbb{N}_0.\eqno(3.1)$$

Indeed, basing on differential-difference equations (2.13), (2.14), we have  

{\bf Theorem 4.} {\it Let $\alpha > -1,\ \nu > 0, \  \lambda, t,  \ge 0, \ \lambda^2+t^2 \neq 0, n \in \mathbb{N}$.  Orthogonal polynomials $P_{n} (x; \lambda, t) $ satisfy the  following integral-difference equations}

$$ P_{n} (x; \lambda, t) = \int_0^\lambda \exp \left( {1\over 2} \int_\xi^\lambda B_n(y,t) dy \right) A_n(\xi,t) P_{n-1}(x; \xi, t) d\xi + \exp \left( {1\over 2} \int_0^\lambda B_n(y,t) dy \right) P_{n} (x; 0, t),\eqno(3.2)$$

$$ P_{n} (x; \lambda, t) = \int_0^t  \left[ x {\partial\over \partial x }  [ P_{n} (x; \lambda, y) ]  -  n P_{n} (x; \lambda, y) -\lambda A_n(\lambda, y)  P_{n-1} (x; \lambda, y)\right] {dy\over  y a_n(\lambda, y)} $$

$$  +  (-1)^n  \lambda^{-n- (1+\alpha)/2}   \bigg[ n! \Gamma (n+\alpha+1) \Gamma(\nu)\bigg]^{1/2} a_n(\lambda, t) \hat{L}^{\alpha}_{n} (x; \lambda).\eqno(3.3)$$

\begin{proof}  The first equation (3.2) is obtained, solving the non-homogeneous ordinary differential-difference equation (2.13) in terms of $\lambda$ with respect to  $P_n$,  and using identity (2.20).  Concerning the second equation, we will resolve (2.14)  under the initial condition (cf. (2.38)) $P_{n} (x; \lambda, 0)= \hat{L}^{\alpha}_{n} (x; \lambda)$. Therefore it is straightforward to write from (2.14) via integration the equality

$$ P_{n} (x; \lambda, t) = \int_0^t  \left[ x {\partial\over \partial x }  [ P_{n} (x; \lambda, y) ]  -  n P_{n} (x; \lambda, y) -\lambda A_n(\lambda, y)  P_{n-1} (x; \lambda, y)\right] {dy\over  y a_n(\lambda, y)} $$

$$  +   {a_n(\lambda, t)\over a_n(\lambda, 0)} \hat{L}^{\alpha}_{n} (x; \lambda),\eqno(3.4)$$
where the integral converges via (2.14) and since $a_n\neq 0.$ Our goal now is to find the value $a_n(\lambda,0)$.  At the same time the three term recurrence relation (1.16) for the modified Laguerre polynomials can be written explicitly,  because due to (2.21), (2.23) and the orthogonality for the normalized Laguerre polynomials we find

$$ B_n(\lambda, 0)= {1\over \lambda} [ 2n+\alpha+1],\quad \lambda > 0,\eqno(3.5)$$

$$A_n(\lambda,0) = {A_n(1,0)\over \lambda}= -  {(n(n+\alpha))^{1/2}\over \lambda},\quad \lambda > 0.\eqno(3.6)$$
Hence, taking into account (1.1) with $t=0$,  (1.17) and   (3.6), we obtain

$$ \prod_{k=1}^n A_k (\lambda, 0)  = {a_0(\lambda,0)\over a_n(\lambda,0) } = {(-1)^n \over \lambda^n } (n! (1+\alpha)_n)^{1/2},\  \lambda >0,\eqno(3.7)$$

$$ a_0(\lambda,0) = {  \lambda^{(1+\alpha)/2} \over (\Gamma(\nu) \Gamma(1+\alpha))^{1/2}},\quad \nu, \lambda  > 0,\eqno(3.8)$$
where $(z)_n$ is the Pochhammer symbol.  Thus it yields the formula

$$  a_n(\lambda,0)  =  { (-1)^n \lambda^{n+ (1+\alpha)/2}  \over \left(n! \Gamma (n+\alpha+1) \Gamma(\nu)\right)^{1/2} }.\eqno(3.9)$$
Substituting this value in (3.4), we end up with (3.3), completing the proof of Theorem 4.

\end{proof}

Concerning the composition orthogonality of the sequence  $\{ P_{n} \}_{n\ge 0}$, we establish the  following result.

{\bf Theorem 5}. {\it Let $\alpha > -1,\ \lambda, t, \nu  \ge 0, \ \lambda^2+t^2 \neq 0.$ The sequence  $\{ P_{n} \}_{n\ge 0}$ is compositionally orthogonal in the sense of Laguerre, i.e.  the corresponding orthogonality conditions have  the form 

$$\int_0^\infty  y^{\nu} e^{-y }  P_n \left({\theta\over t}; \lambda, t \right)  \theta^m \left\{ {\Gamma(1+\alpha)\ y^\alpha \over (\lambda y+ t)^{\alpha +1}} \right \} dy = 0,\quad m=0,1,\dots, n-1,\ \quad n \in \mathbb{N},\eqno(3.10)$$
where $\theta = y D y$.}

\begin{proof}  Recalling integral representation (1.10) of the function $\rho_\nu$ in terms of Laguerre polynomials, we substitute its right-hand side  in (1.2) and interchange the order of integration by Fubini's theorem via the absolute convergence of the iterated integral.  Hence we write

$$\int_{0}^{\infty}  P_{n} (x; \lambda, t)  e^{-\lambda x} \rho_\nu(x t) \  x^{\alpha+m} dx =  (-1)^m m! \int_0^\infty y^{\nu+m -1} e^{-y }  L_m^\nu(y) \int_{0}^{\infty}  P_{n} (x; \lambda, t)  e^{- x (\lambda + t/y) } x^{\alpha} dx dy.\eqno(3.11)$$
Then, invoking the Rodrigues formula form Laguerre polynomials, we observe (see (3.1))

$$ m!  y^{\nu+m} e^{-y }  L_m^\nu(y) = \theta^m \left\{ y^{\nu} e^{-y }\right\}.\eqno(3.12)$$
Moreover, the differentiation by $y$ under the integral sign via the absolute and uniform convergence allows to deduce the following operator formula

$$ {1\over y} \int_{0}^{\infty}  P_{n} (x; \lambda, t)  e^{- x (\lambda + t/y) } x^{\alpha} dx = P^\nu_n \left({\theta\over t}; \lambda, t \right)  \left\{ {1\over y}  \int_{0}^{\infty}  e^{-x(\lambda +t/y)} x^\alpha dx\right \}$$

$$=   P^\nu_n \left({\theta\over t}; \lambda, t \right)  \left\{ {\Gamma(1+\alpha)\ y^\alpha \over (\lambda y+ t)^{\alpha +1}} \right \}.\eqno(3.13)$$
Plugging (3.12), (3.13) in the right-hand side of (3.11), we integrate by parts, eliminating the integrated terms. Thus we derive

$$ \int_{0}^{\infty}  P_{n} (x; \lambda, t)  e^{-\lambda x} \rho_\nu(x t) \  x^{\alpha+m} dx = (-1)^m \int_0^\infty \theta^m \left\{ y^{\nu} e^{-y }\right\}  P^\nu_n \left({\theta\over t}; \lambda, t \right)  \left\{ {\Gamma(1+\alpha)\ y^\alpha \over (\lambda y+ t)^{\alpha +1}} \right \} dy$$

$$=  \int_0^\infty  y^{\nu} e^{-y }  P^\nu_n \left({\theta\over t}; \lambda, t \right)  \theta^m \left\{ {\Gamma(1+\alpha)\ y^\alpha \over (\lambda y+ t)^{\alpha +1}} \right \} dy.$$
This implies (3.10) and completes the proof.

\end{proof}

\section{One-parameter case $\lambda = 1-t$}

Let us consider the sequence of orthogonal polynomials $\{P_{n}(x; 1-t, t)\}_{n \ge 0}$ with respect to the weight $e^{- (1-t)x} \rho_\nu(x t),\ 0 \le t \le 1$.  The respective orthogonality conditions (1.2) will take the form 

$$\int_{0}^{\infty}  P_{n} (x; 1-t, t)  e^{- (1-t) x} \rho_\nu(x t) \  x^{\alpha+m} dx = 0,\quad m=0,1,\dots, n-1,\ \quad n \in \mathbb{N}.\eqno(4.1)$$
Then, making a differentiation with respect to parameter,  we deduce from (4.1) similar to (2.5), (2.6)

$$\int_{0}^{\infty}  {\partial\over \partial t } [ P_{n} (x; 1-t, t) ] e^{-(1-t) x} \rho_\nu(x t) \  x^{\alpha+m} dx +  \int_{0}^{\infty}   P_{n} (x; 1-t, t) e^{-(1-t) x} \rho_\nu(x t) \  x^{\alpha+m+1} dx$$

$$- \int_{0}^{\infty}  P_{n} (x; 1-t, t)  e^{- (1-t) x} \rho_{\nu-1} (x t) \  x^{\alpha+m+1} dx = 0,\quad m=0,1,\dots, n-1.\eqno(4.2)$$
Now,   employing again (1.15) and orthogonality (4.1), we find from the latter equalities (4.2)

$$\int_{0}^{\infty}   t {\partial\over \partial t } [ P_{n} (x; 1-t, t) ] e^{-(1-t) x} \rho_\nu(x t) \  x^{\alpha+m} dx + t  \int_{0}^{\infty}   P_{n} (x; 1-t, t) e^{-(1-t) x} \rho_\nu(x t) \  x^{\alpha+m+1} dx$$

$$- \int_{0}^{\infty}  P_{n} (x; 1-t, t)  e^{- (1-t) x} \rho_{\nu+1} (x t) \  x^{\alpha+m} dx = 0,\quad m=0,1,\dots, n-1.\eqno(4.3)$$
On the other hand, integrating by parts in (4.1) and using (1.14), (1.15) and the orthogonality, it gives (cf. (2.7))

$$ (1-t) \int_{0}^{\infty}   P_{n} (x; 1-t, t) e^{- (1-t) x} \rho_\nu(x t) \  x^{\alpha+m+1} dx +  \int_{0}^{\infty}   P_{n} (x; 1-t, t) e^{-(1-t) x} \rho_{\nu+1} (x t) \  x^{\alpha+m} dx $$

$$- \int_{0}^{\infty}  {\partial\over \partial x } [ P_{n} (x; 1-t, t) ] e^{- (1-t) x} \rho_\nu(x t) \  x^{\alpha+m+1} dx = 0,\quad m=0,1,\dots, n-1.\eqno(4.4)$$
Eliminating the integral with $\rho_{\nu+1}$ via (4.3),  we get

$$ \int_{0}^{\infty}  \left( t {\partial\over \partial t} -  x {\partial\over \partial x } \right)  [ P_{n} (x; 1-t, t) ] e^{- (1-t) x} \rho_\nu(x t) \  x^{\alpha+m} dx$$

$$+   \int_{0}^{\infty}   P_{n} (x; 1-t, t) e^{- (1-t) x} \rho_\nu(x t) \  x^{\alpha+m+1} dx = 0,\quad m=0,1,\dots, n-1.\eqno(4.5)$$
This yields immediately an analog of the differential-difference equation (2.14)

$$ \left( t {\partial\over \partial t} -  x {\partial\over \partial x } \right)  [ P_{n} (x; 1-t, t) ]  - \bigg[ { t \over a_n}\  {\partial a_n \over \partial t } -  n \bigg] P_{n} (x; 1-t, t) + A_n P_{n-1} (x; 1-t, t) = 0.\eqno(4.6)$$
Equating coefficients in front of $x^{n-1}$ in (4.6), we derive the Toda-type equations (cf. (2.18), (2.19))

$$   {\partial  \over \partial t }  \left( { t b_n \over a_n} \right)  = -  A_n^2,\eqno(4.7)$$

$$  {\partial  \over \partial t }  \left( t B_n\right)  = \left[ A_{n+1}^2- A_{n}^2\right].\eqno(4.8)$$
Finally, as in (2.39), (2.40) we conclude from (4.5) that the sequence $\left\{\left( t {\partial\over \partial t} -  x {\partial\over \partial x } \right)  [ P_{n} (x; 1-t, t) ] \right\}_{n\ge 0}$ is quasi-orthogonal with respect to the weight $e^{-(1-t) x} \rho_\nu(x t) x^{\alpha}$.

\bibliographystyle{amsplain}

\end{document}